\documentclass[11pt, reqno]{amsart} 
\usepackage{amsmath, amsthm, amsfonts, amssymb} 
 \usepackage[colorlinks]{hyperref} 
 \usepackage[T1]{fontenc} 
 \usepackage[mathscr]{eucal} 
 \usepackage[polish,english]{babel} 
 \newtheorem{theorem}{Theorem}[section]

 \newtheorem{corollary}[theorem]{Corollary} 
  
 \newtheorem{definition}[theorem]{Definition} 
 \newtheorem{example}[theorem]{Example} 
  
 \newtheorem{lemma}[theorem]{Lemma}

 \newtheorem{proposition}[theorem]{Proposition}

 \usepackage{amsmath, amsthm, amscd, amsfonts, amssymb, graphicx, color}

\usepackage{lipsum}
 
\begin{document} \title[EDGE QUASI $\lambda$-DISTANCE-BALANCED GRAPHS IN METRIC SPACE] {EDGE QUASI $\lambda$-DISTANCE-BALANCED GRAPHS IN METRIC SPACE} 
 \date{} \author[Zohreh Aliannejadi and Somayeh Shafiee Alamoti ]{Zohreh Aliannejadi and Somayeh Shafiee Alamoti} 
 \begin{abstract} 
In a graph $A$, the measure $|M_g^A(f)|=m_g^A(f)$ for each arbitrary edge $f=gh$ counts the edges in $A$ closer to $g$ than $h$. $A$ is termed an edge quasi-$\lambda$-distance-balanced graph in a metric space (abbreviated as $EQDBG$), where a rational number ($>1$) is assigned to each edge $f=gh$ such that $m_g^A(f)=\lambda^{\pm1}m_h^A(f)$. This paper introduces and discusses these graph concepts, providing essential examples and construction methods. The study examines how every $EQDBG$ is a bipartite graph and calculates the edge-Szeged index for such graphs. Additionally, it explores their properties in Cartesian and lexicographic products. Lastly, the concept is extended to nicely edge distance-balanced and strongly edge distance-balanced graphs revealing significant outcomes.
 \end{abstract} 
 \maketitle \textbf{Mathematics Subject Classification (2010): 05C12; 05C25.}
 \section{Introduction}
The concept of graphs plays a crucial role in modeling various phenomena and has been extensively utilized in numerous studies in recent years. Graph theory is particularly valuable in categorizing graphs based on their distinguishing characteristics. This caregorization is prominently highlighted in distance-balanced graphs, as discussed in [13]. Additionally, this topic has been explored in several papers, and for further detailes, we recommend referring to ([3],[4],[6],[11],[14]-[19]) and the related references.\\ 
Throughout this paper, we focus on a connected, finite, and undirected graph $A$, where the vertex set is denoted as $V(A)$ and the edge set as $E(A)$. $A$ is referred to as a metric space if it is satisfies the triangle inequality: $d(g,i) \leqslant d(g,h) + d(h,i)$ for every vertex $g, h, i \in V(A)$. A graph $A=(A,d)$ is considered a metric space if there exists a metric space $M=(A, \acute{d})$ such that $d(g,h)=\acute{d}(g,h)$ for every vertex $g, h$ in $V(A)$. In graph $A$, the distance between vertices $g$ and $h$ in $V(A)$ is defined as the number of edges in the shortest path connecting them, denoted as $d_A(g,h)$. For any two given vertices $g$ and $h$ in $V(A)$, we define $n_g^A(f)=|W_{a,h}^A|=|\{g\in V(A)| d_A(a,g)< d_A(a,h)\}|$. Similarly, $n_h^A(f)=|W_{h,g}^A|$, where $A$ is termed quasi-$\lambda$-distance-balanced (denoted as $\lambda-QDB$). For adjacent vertices $g$ and $h$ in $A$, a constant 
$\lambda>1$ is set, such that $|W_{g,h}^A|=\lambda^{ \pm 1} |W_{h,g}^A|$.\\
In this paper, we aim to extend the previous definition of edges in graph $A$. Initially, we establish the definition of the distance between two edges in a metric space. Subsequently, we introduce the distance between $f$ and $\acute{f}$ for any two desired edges $f=gh$ and $\acute{f}=\acute{g}\acute{h}$. 
\begin{center}
$d_A(f,\acute{f})=\min\{d_A(g,\acute{f}), d_A(h,\acute{f})\}=\min\{d_A(g,\acute{g}),d_A(g,\acute{h}),d_A(h,\acute{g}),d_A(h,\acute{h})\}$.
\end{center}
Set $\hspace*{0.1cm}$ $M_g(f)=\{\acute{f}\in E(A)| d_A(g,\acute{f})< d_A(h,\acute{f})\}$ and $m_g(f)=|M_g(f)|$,\\
$\hspace*{0.7cm}$ $M_h(f)=\{\acute{f}\in E(A)| d_A(h,\acute{f})< d_A(g,\acute{f})\}$ and $m_h(f)=|M_h(f)|$,\\
and $\hspace*{0.0001cm}$ $M_0(f)=\{\acute{f}\in E(A)| d_A(g,\acute{f})= d_A(h,\acute{f})\}$ and $m_0(f)=|M_0(f)|$.\\
 
Presume that $f=gh\in E(A)$. For every two integers $i,j$ we consider:
\begin{center}
$\acute{D}_j^i(f)=\{\acute{f}\in E(A)|d_A(\acute{f},g)=i, d_A(\acute{f},h)=j\}$.
\end{center}
The sets $\acute{D}_i^i(f)$ are used to form a "distance partition" of $E(A)$ for the edge $f=gh$. Only the sets $\acute{D}_i^{i-1}(f)$, $\acute{D}_i^i(f)$ and, $\acute{D}_{i-1}^i(f)$, where $(1\leqslant i\leqslant d)$ may be non-empty based on the triangle inequality (where $d$ is the diameter of the graph $A$), while $\acute{D}_0^0(f)$ is an empty set. It should be noted that the sets $W_{g,h}^A$ and $M_g^A(f)$ are known as the Szeged index and edge-Szeged index of graph $A$ in chemical graph theory, respectively. These are defined as $S_z(A)=\sum_{gh\in E(A)}|W_{g,h}|.|W_{h,g}|$ and $S_{{z}_{e}}(A)=\sum_{gh\in E(A)}m_{g}^A(f).m_{h}^A(f)$, respectively, as described in [7,8,9].\\

A graph $A$ is considered to be edge distance-balanced ($EDB$) if $m_g^A(f)=m_h^A(f)$. Additionally, $A$ is referred to as nicely edge distance-balanced ($NEDB$) if there exists a positive integer $A$ such that $m_g^A(f)=m_h^A(f)=\acute{\gamma}_A$ for every edge $f=gh$. A graph is classified as strongly edge distance-balanced ($SEDB$) if for each $i\geq 1$, $\acute{D}_{i-1}^i(f)$ is equal to $\acute{D}_i^{i-1}(f)$. It is evident from the given definition that if the graph is $SEDB$, then it is an $EDB$ graph.\\

A graph denoted as $A$ is considered edge quasi-distance-balanced in a metric space, where $\lambda$, a positive rational number exceeding $1$, satisfies the condition that for every edge $f=gh$ in $E(A)$, either $m_g(f)=\lambda m_h(f)$ or $m_h(f)=\lambda m_g(f)$. To simplify, we use $\lambda-EQDBG$ in place of edge quasi-$\lambda$-distance-balanced graph and $\lambda-EQDBG$ for short.\\
The theorem asserts that there does not exist a non-bipartite $\lambda-EQDBG$.
\begin{theorem}\label{1.1}
Presume that $A$ is a $\lambda-EQDBG$. Then $A$ is a bipartite graph $K_{m,n}$, where $m\neq n$.
\end{theorem}

 \begin{theorem}\label{1.2}
 A bipartite graph $A$ is $\lambda-EQDB$ if and only if
 \begin{center}
 $S_{{z}_{e}}(A)=\lambda.\frac{(\lambda n+\lambda+1)^2.\lambda n^2}{(\lambda+1)^2}$.
 \end{center}
 \end{theorem}

We characterize situations, in which the lexicographic and catresian products would lead to a $\lambda-EQDBG$ in the below.

\begin{theorem}\label{1.3}
Let $A$ and $B$ be graphs. Then the catresian product $A\square B$ is $\lambda-EQDB$ if and only if both $A$ and $B$ are both $\lambda-EQDB$ and $\lambda-QDB$.
\end{theorem}

\begin{theorem}\label{1.4}
Assume that $A$ and $B$ are graphs. Then the lexicographic product $A[B]$ is $\lambda-EQDB$ if and only if $A$ is $\lambda-EQDB$ and also $B$ is an empty graph.
\end{theorem}

We sort out this article as follows. In the next two sections, we reveal some related facts about being $\lambda-EQDB$ and our main results regarding characteristics of $\lambda-EQDB$ graphs amidst the product graphs are studied. In Section 4 we are going to introduce some methods to construct concerning $\lambda-EQDB$ graphs. Finally, we define two new graphs nicely edge quasi $\lambda$-distance-balanced and strongly edge quasi $\lambda$-distance-balanced and state intersting results.

 \section{Main Results} 
 In this segment, the principal result already explained in the introduction will be proved. We would express proof of Theorem 1.1.\\\\
\textbf{Proof of Theorem 1.1} Inspired of proof of [2. Theorem 1], presume that $A$ is $\lambda-EQDB$ such that $d=diam(A)$, and the edge set $\{f_1,f_2,...,f_{2m+1}\}$ organize an odd cycle that its length is $2m+1$ for $f=gh\in E(A)$ and
\begin{center}
$P_{ij}=\{f\in E(A)|d(f,f_{i+k})=s_{jk},$\hspace*{0.3cm}$ s_{jk}=\{1,2,...,d\}, k=0,1,...,2m\}$,$\hspace*{0.3cm}$ $2\leqslant j\leqslant n$,
\end{center}
so that $M_g(f_i)=(\bigcup_{j=1}^nP_{ij})\cup\{f_{i+2m}\}$ and
$M_h(f_i)=(\bigcup_{j=1}^nP_{(i+1)j})\cup\{f_{i+2}\}$, in which mod $2m+1$ are effectuated by the computations indexes $i$ for $n\in\mathbb{N}$. Let $|P_{ij}|=p_{ij}$ for $i=0,1,...,2m$ and $j=1,2,...,n$ and the below assumption there is $t_i\in \{ \pm1\}$, $i=0,1,...,2m$, therefore,
\begin{center}
$\sum_{j=1}^n p_{0j}+1= \lambda^{t_{0}}(\sum_{j=1}^n p_{1j}+1)$,
\end{center}
\begin{center}
$\sum_{j=1}^n p_{1j}+1= \lambda^{t_{1}}(\sum_{j=1}^n p_{2j}+1)$,
\end{center}
$\hspace*{5.70cm}$ .\\
$\hspace*{5.70cm}$ .\\
$\hspace*{5.70cm}$ .\\
\begin{center}
$\hspace*{0.4 cm}$ $\sum_{j=1}^n p_{(2m-1)j}+1= \lambda^{t_{2m-1}}(\sum_{j=1}^n p_{(2m)j}+1)$,
\end{center}

\begin{center}
$\sum_{j=1}^n p_{(2m)j}+1= \lambda^{t_{2m}}(\sum_{j=1}^n p_{0j}+1)$.
\end{center}
Now, it is concluded from the combination of all the $(2m+1)$ above formulas that $\lambda^{\Sigma_{i=0}^{2m}t_i}=1$, that is,
$\Sigma_{i=0}^{2m}t_i=0$. In other words,
 \begin{center}
 $t_i\in \{\pm1\}\Rightarrow 1\leqslant|\Sigma_{i=0}^{2m}t_i|$.
 \end{center}
 Then, it is a contradiction and so there is no odd cycle in $A$. The proof is completed.\hfill\qed
\begin{corollary}\label{2.1}
For positive integer quantities $m<n$, complete bipartite graph $K_{m,n}$ is a $\left(\frac{n-1}{m-1}\right)$-$EQDBG$.
\end{corollary}
\proof Assume that $X=\{x_1,x_2,...,x_m\}$ and $Y=\{y_1,y_2,...,y_n\}$ are two sets of complete bipartite graph $K_{m,n}$ and $f_{ij}=x_iy_j$ is an arbitrary edge, for each $i\in\{1,2,...,m\}$ and for each $j\in\{1,2,...,n\}$. Then we see vividly
\begin{center}
	$M_{x_i}(f_{ij})=M_{y_i}(f_{ij})=\{f_{11},...,f_{1n},f_{2n},...,f_{m1},...,f_{mn}\}\backslash \{f_{ij}\}$.
\end{center}
So, $m_{x_i}(f_{ij})=n-1$ and $m_{y_i}(f_{ij})=m-1$. We conclude that $K_{m,n}$ is a $\left(\frac{n-1}{m-1}\right)$-$EQDBG$.\hfill\qed\\\\
\textbf{Proof of Theorem 1.2} Presume that $A$ is a $\lambda-EQDBG$. For every edge $f=gh$, we have $m_g^A(f)=\lambda m_h^A(f)$ and since $A$ is bipartite also $m_g^A(f)+m_h^A(f)+m_0^A(f)=|E(A)|-1$ holds. Thus,
 \begin{center}
$m_g^A(f)+m_h^A(f)=|E(A)|-m_0^A(f)-1$.
\end{center}
Then,
\begin{center}
$\lambda^2(m_h^A(f))^2+(m_h^A(f))^2+2\lambda(m_h^A(f))^2=(|E(A)|-m_0^A(f)-1)^2$.
\end{center}
Hence,
\begin{center}
$(m_h^A(f))^2=\frac{(|E(A)|-m_0^A(f)-1)^2}{(\lambda+1)^2}$,
\end{center}
and so,
 \begin{center}
 $\sum_{f\in E(A)}m_{g}^A(f).m_{h}^A(f)=\lambda\sum_{f\in E(A)}(m_{h}^A(f))^2=$
 \end{center}
 \begin{center}
 $\lambda.\frac{(|E(A)|-m_0^A(f)-1)^2.|E(A)|}{(\lambda+1)^2}=\lambda.\frac{(\lambda n+\lambda+1)^2.\lambda n^2}{(\lambda+1)^2}$.
 \end{center}
For converse, let $S_{{z}_{e}}(A)=\lambda.\frac{(\lambda n+\lambda+1)^2.\lambda n^2}{(\lambda+1)^2}$. Since $A$ is a graph, we have 
$m_g^A(f)+m_h^A(f)+m_0^A(f)=|E(A)|-1$ and also 
\begin{center}
$m_0^A(f)+(\lambda+1)m_h^A(f)=|E(A)|-1$.
\end{center}
Hence,
\begin{center}
$m_h^A(f)=\frac{(|E(A)|-m_0^A(f)-1)}{(\lambda+1)}$.
\end{center}
As well as 
\begin{center}
$m_g^A(f).m_h^A(f)=\lambda.\frac{(|E(A)|-m_0^A(f)-1)^2}{(\lambda+1)^2}$,
\end{center}
and hence
\begin{center}
$m_g^A(f)=\lambda m_h^A(f)=\lambda.\frac{(|E(A)|-m_0^A(f)-1)}{(\lambda+1)}$.
\end{center}
This implies that $A$ is a $\lambda-EQDBG$.\hfill\qed\\

\section{ $\lambda-EQDB$ Property in product graphs}

In examining scenarios where the \emph{Cartesian product} results in a $\lambda-EQDBG$, we note that such product graphs, formed by graphs $A$ and $B$, have a vertex set denoted as $V(A\square B)=V(A)\times V(B)$. Let $(a_1,b_1)$ and $(a_2,b_2)$ represent distinct vertices in $V(A\square B)$. In the Cartesian product $A\square B$, if vertices $(a_1,b_1)$ and $(a_2,b_2)$ coincide in one coordinate and are adjacent in the other coordinate, then they are considered adjacent. Specifically, either $a_1=a_2$ and $b_1b_2\in E(B)$, or $b_1=b_2$ and $a_1a_2\in E(A)$. Clearly, for vertices, we have: 
\begin{center}
$d_{A\square B}((a_1,b_1),(a_2,b_2))=d_A(a_1,a_2)+d_B(b_1,b_2)$.
\end{center}
For edges we have:
\begin{center}
$d_{A\square B}((a,b)(a_1,b_1),(\acute{a},\acute{b})(\acute{a}_1,\acute{b}_1))=$
\end{center}
\begin{center}
$\min\{d_{A\square B}((a,b),(\acute{a},\acute{b})),d_{A\square B}((a,b),(\acute{a}_1,\acute{b}_1)),d_{A\square B}((a_1,b_1),(\acute{a},\acute{b})$, \\ $d_{A\square B}((a_1,b_1),(\acute{a}_1,\acute{b}_1))\}=$
 \end{center}
 \begin{center}
$\min\{d_A(a,\acute{a})+d_B(b,\acute{b}),d_A(a,\acute{a}_1)+d_B(b,\acute{b}_1),
d_A(a_1,\acute{a})+d_B(b_1,\acute{b}),d_A(a_1,\acute{a}_1)+d_B(b_1,\acute{b}_1)\}$.
\end{center}

We now verify the $\lambda-EQDB$ property of cartesian product graphs. We demonstrate Theorem 1.3.\\\\
\textbf{Proof of Theorem 1.3} Let $a_1,a_2$ be adjacent vertices in $V(A)$, and $b_1,b_2$ are two adjacent vertices in $V(B)$. We assume that $(a_1,b_1), (a_2,b_1),\\ (a_1,b_2)\in V(A\square B)$. We observe that
\begin{center}
$(a,b)(\acute{a},\acute{b})\in M_{(a_1,b_1)}((a_1,b_1)(a_2,b_1))\Leftrightarrow$
\end{center}
\begin{center}
$\min\{d_{A\square B}((a,b),(a_1,b_1)),d_{A\square B}((\acute{a},\acute{b}),(a_1,b_1))\}<\min\{d_{A\square B}((a,b),(a_2,b_1)),d_{A\square B}((\acute{a},\acute{b}),(a_2,b_1))\}\Leftrightarrow$
\end{center}
\begin{center}
$\min \{d_A(a,a_1)+d_B(b,b_1),d_A(\acute{a},a_1)+d_B(\acute{b},b_1)\}<\min \{d_A(a,a_2)+d_B(b,b_1),d_A(\acute{a},a_2)+d_B(\acute{b},b_1)\}\Leftrightarrow$
\end{center}
\begin{center}
$\min\{d_A(a,a_1),d_A(\acute{a},a_1)\}<\min\{d_A(a,a_2),d_A(\acute{a},a_2)\}$.
\end{center}
We deduce that
\begin{center}
$M_{(a_1,b_1)}((a_1,b_1)(a_2,b_2))=\{(a,b)(\acute{a},\acute{b})\in E(A\square B)|a\acute{a}\in E(A) , b=\acute{b}$ or\\ $b\acute{b}\in E(B)$ , $a=\acute{a}$ ,
\end{center}

$\hspace*{0.5cm}$ $\min\{d_A(a,a_1),d_A(\acute{a},a_1)\}<\min\{d_A(a,a_2),d_A(\acute{a},a_2)\}\}$. $\hspace*{1.60cm}$(3.1)\\\\
The application of (1) and [19, Theorem 2.1] leads to the following conclusions:\\\\
$\hspace*{0.5cm}$ $m^{A\square B}_{(a_1,b_1)}((a_1,b_1)(a_2,b_2))=m^A_{a_1}(a_1a_2).|V(B)|+n^A_{a_1}(a_1a_2).|E(B)|$.$\hspace*{0.65cm}$(3.2)\\\\
Similarly\\\\
$\hspace*{0.5cm}$ $m^{A\square B}_{(a_2,b_2)}((a_2,b_2)(a_1,b_1))=m^A_{a_2}(a_1a_2).|V(B)|+n^A_{a_2}(a_1a_2).|E(B)|$,$\hspace*{0.4cm}$(3.3)\\\\
$\hspace*{0.6cm}$ $m^{A\square B}_{(a_1,b_1)}((a_1,b_1)(a_1,b_2))=m^B_{b_1}(b_1b_2).|V(A)|+n^B_{b_1}(b_1b_2).|E(A)|$,$\hspace*{0.6cm}$(3.4)\\\\
$\hspace*{0.6cm}$ $m^{A\square B}_{(a_1,b_2)}((a_1,b_2)(a_1,b_1))=m^B_{b_2}(b_1b_2).|V(A)|+n^B_{b_2}(b_1b_2).|E(A)|$.$\hspace*{0.6cm}$(3.5)\\\\
If both $A$ and $B$ are both $\lambda-EQDBG$ and $\lambda-QDBG$, then by comparing (3.2) and (3.3), it is derived that:\\\\
$m^A_{a_1}(a_1a_2).|V(B)|+n^A_{a_1}(a_1a_2).|E(B)|=\mu m^A_{a_2}(a_1a_2).|V(B)|+\mu n^A_{a_2}(a_1a_2).|E(B)|$,\\\\
for $\mu \in\{\lambda,\frac{1}{\lambda}\}$, therefore
\begin{center}
$m^{A\square B}_{(a_1,b_1)}((a_1,b_1)(a_2,b_2))=\mu m^{A\square B}_{(a_2,b_2)}((a_1,b_1)(a_2,b_2))$.
\end{center}
By analogy, using (3.4) and (3.5), we conclude that
\begin{center}
$m^{A\square B}_{(a_1,b_1)}((a_1,b_1)(a_1,b_2))=\mu m^{A\square B}_{(a_1,b_2)}((a_1,b_1)(a_1,b_2))$,
\end{center}
and hence $A\square B$ is $\lambda-EQDB$.\\
For converse, assume that $A\square B$ is $\lambda-EQDB$. Then applying (3.2) and (3.3) we attain
\begin{center}
$m^{A\square B}_{(a_1,b_1)}((a_1,b_1)(a_2,b_2))=\mu m^{A\square B}_{(a_2,b_2)}((a_1,b_1)(a_2,b_2))\Rightarrow$
\end{center}
$m^A_{a_1}(a_1a_2).|V(B)|+n^A_{a_1}(a_1a_2).|E(B)|=\mu m^A_{a_2}(a_1a_2).|V(B)|+\mu n^A_{a_2}(a_1a_2).|E(B)|$,\\\\
for $\mu \in\{\lambda,\frac{1}{\lambda}\}$. Hence $A$ is $\lambda-EQDB$ and $\lambda-QDB$. 
This implies $A$ is $\lambda-EQDB$ and $\lambda-QDB$. Similarly, by comparing (3.4) and (3.5), it's concluded that $B$ is $\lambda-EQDB$ and $\lambda-QDB$. This completes the result.\hfill\qed\\

The \emph{lexicographic product} graphs, denoted as $A[B]$, are defined for the graphs $A$ and $B$. The vertex set of $A[B]$ is given by $V(A[B])=V(A)\times V(B)$, and two vertices $(a_1,b_1)$ and $(a_2,b_2)$ are considered adjacent in the lexicographic product if $a_1a_2\in E(A)$ or if $a_1=a_2$ and $b_1b_2\in E(B)$ (for further details, refer to [12, p. 22]). Since $A$ is a graph, it is evident for vertices that
\begin{center}
$\ d_{A[B]}((a_1,b_1),(a_2,b_2))=\left\{ \begin{array} {cc} d_A(a_1,a_2) & if a_1\neq a_2\\ \min\{2,d_B(b_1,b_2)\} & if a_1=a_2. \end{array}\right.$
\end{center}
And for edges we have:
\begin{center}
$ d_{A[B]}((a,b)(a_1,b_1),(\acute{a},\acute{b})(\acute{a}_1,\acute{b}_1))=$
\end{center}

\begin{equation}\min 
\left\{
\begin{array}{cccc}
d_A(a,\acute{a}) & if a\neq \acute{a}, & \min\{2,d_B(b,\acute{b})\} & if a=\acute{a} \notag\\
d_A(a,\acute{a_1}) & if a\neq \acute{a_1}, & \min\{2,d_B(b,\acute{b_1})\} & if a=\acute{a_1}\\
d_A(a_1,\acute{a}) & if a_1\neq \acute{a}, & \min\{2,d_B(b_1,\acute{b})\} & if a_1=\acute{a}\\
d_A(a_1,\acute{a}_1) & if a_1\neq \acute{a}_1, & \min\{2,d_B(b_1,\acute{b}_1)\} & if a_1=\acute{a}_1
\end{array}
\right\}.
\end{equation}

It is demonstrated that the lexicographic product A[B] of graphs $A[B]$ of graphs $A$ and $B$ is $DB$ if and only if $A$ is $DB$ and $B$ is regular inspired by [13, Theorem 4.2]. This notion is extended to $\lambda-EQDB$ graphs in the following theorem. Let's proceed to the proof of Theorem 1.4.\\\\
Let's consider $A$ as a connected graph, and $A[B]$ as a $\lambda-EQDBG$. By Theorem 1, it follows that $A[B]$ is a bipartite graph. Since $B$ has at least one edge, it becomes evident that $A[B]$ is not a bipartite graph. Therefore, $B$ must be an empty graph. Assume $a_1$ and $a_2$ are adjacent vertices in $A$, and $b_1,b_2\in V(B)$ ($b_1$  and $b_2$ are not necessarily distinct). According to the definition of the lexicographic product, $(a_1,b_1)$ and $(a_2,b_2)$ are adjacent and we observe that,
\begin{center}
$M^{A[B]}_{(a_1,b_1)}=\{(a_1,b_1)(a,b)\}\cup \{(a_2,\acute{b})(a,b)|\acute{b}\in  V(B)\backslash \{b_2\}\}$
\end{center}
\begin{center}
$\cup \{(a,b)(a_1,\acute{b})| aa_1\in M^A_{a_1}(a_1a_2), a\acute{a}\neq a_1a_2, \acute{b}\in V(B)$\},
\end{center}
where $a\acute{a}\neq a_1a_2$, that is $\{a\neq a_1, \acute{a}\neq a_2\}$ or $\{\acute{a}\neq a_2, \acute{a}\neq a_1\}$.\\
Thus, $m^{A[B]}_{(a_1,b_1)}((a_1,b_1)(a_2,b_2))=|V(B)| . m^A_{a_1}(a_1a_2)$. Similarly we attain\\$m^{A[B]}_{(a_2,b_2)}((a_1,b_1)(a_2,b_2))=|V(B)| . m^A_{a_2}(a_1a_2)$. For an empty graph $B$, it is straightforward to realize that the lexicographic product $A[B]$ is $\lambda-EQDB$ if and only if $A$ is $\lambda-EQDB$. The proof is completed.\hfill\qed

\section{Construction of some $\lambda-EQDB$ graphs}

In this part of paper, we would give some examples of $\lambda-EQDB$ graphs. We already have shown that all $\lambda-EQDB$ graphs are bipartite graphs (Theorem 1.1). The Following examples are constructed from complete bipartite graphs.\\
We need to introduce subdivision-related graph named $S(A)$ and is defined by [5, 20].

\begin{definition}
Presume that $A$ is a connected graph. Then $S(A)$ is the graph attained by placing an extra vertex in each edge of $A$. Equvalently, each edge of $A$ is exchanged by a path of length 2. (For an example, see Fig. 1, 2).
\end{definition}

\begin{figure}
\minipage{0.55\textwidth}
\includegraphics[width=\linewidth]{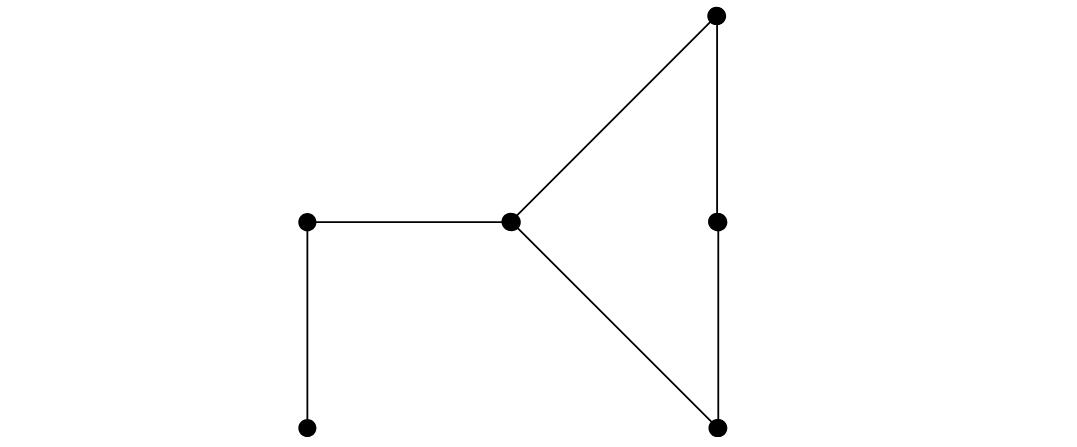}
\caption{Graph $A$}\label{fig:1} 
\endminipage
\minipage{0.55\textwidth}
\includegraphics[width=\linewidth]{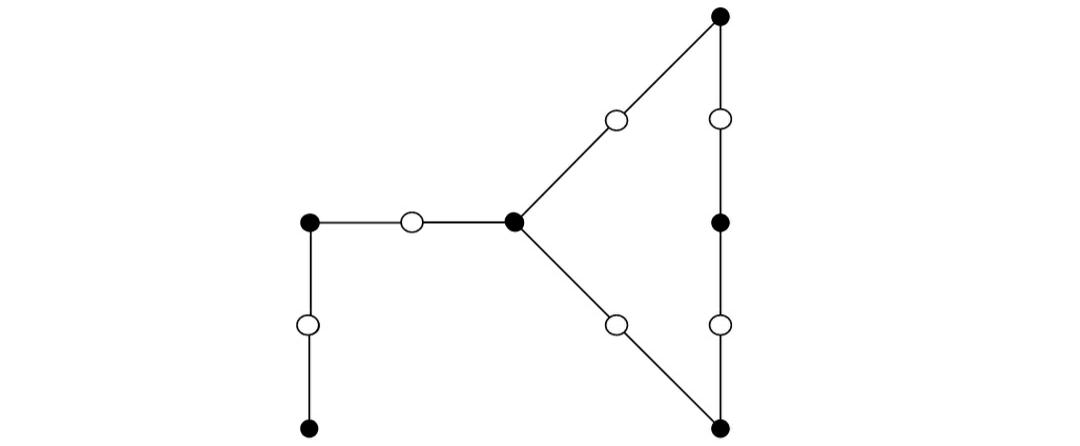}
\caption{Graph $S(A)$}\label{fig:2}
\endminipage
\end{figure}

\begin{proposition}
Presume that $A$ is a regular complete bipartite graph of order $m$, where $m\geqslant 3$. The graph $S(A)$ is $\lambda-EQDB$ if and only if $\lambda=\frac{m^2-1}{2m-1}$.
\end{proposition}

\begin{proof}
We know that $S(A)$ is constructed by two bipartitions sets $X$ and $Z$ of order and size $m$ in $A$ and the set $Y$ of order $m^2$ and size 2 according to the definition of $S(A)$. To prove that graph $S(A)$ is $\lambda-EQDB$, we need to show that $m_x(xy)=\mu m_y(xy)$ and $m_z(yz)=\mu m_y(yz)$ for $x\in X, y\in Y$ and $z\in Z$, and for any pair of edges $xy, yz\in E(S(A))$ and $\mu\in \{\lambda, \frac{1}{\lambda}\}$.
If $A$ is a graph with $m=1$, then $S(A)$ is a path with length 2. So $S(A)$ is not a $\lambda-EQDBG$.\\
It is explicit that if in the graph $A$, $m=2$, then $S(A)$ is an even cycle that is edge distance-balanced. So it can not be $\lambda-EQDB$.\\
We now claim that if $A$ is a graph, where $m=3$ and an edge $f=xy\in E(S(A))$, then it concludes that
\begin{center}
$\{\acute{f} \in S(A)|d(\acute{f}, x) \in \{0,1,2\}\} \in M_x(xy)$,
\end{center}
and
\begin{center}
$\{\acute{f} \in S(A)|d(\acute{f}, y) \in \{0,1,2\}\} \in M_y(xy)$.
\end{center}
Since $d(\acute{f},x)=3$ and $d(\acute{f},y)=3$ for edges $\acute{f}, f\in E(S(A))$, hence $\acute{f}\in M_0(xy)$. It is easily seen that $diam (S(A))=4$. Hence, by the definition of the distance between an edge and a vertex, the maximum distance between them in graph $S(A)$ is 3. For each edges $\acute{f}, f\in E(S(A))$, we investigate the following cases:\\
\textbf{(1)} Let $d(\acute{f},x)=0$. Clearly, the adjacent edges belong to $M_x(xy)$, that is, $m_x(xy)=m-1$.\\
\textbf{(2)} Let $d(\acute{f},x)=1$. Hence, we must compute the number of adjacent edges of $\acute{f}\in M_x(xy)$, with $d(\acute{f},x)=0$. Then, regarding to $deg(y)=2$ for any $y\in Y$, we have $m_x(xy)=m-1$.\\
\textbf{(3)} Now for the last case we consider that $d(\acute{f},x)=2$, that is, the number of adjacent edges of $\acute{f}\in M_x(xy)$, with $d(\acute{f},x)=1$. Therefore, we obtain $m_x(xy)=(m-1)(m-1)$.\\
Then, by the above cases \textbf{(1)}, \textbf{(2)} and \textbf{(3)} for any edge $f=xy\in E(S(A))$ it implies that
\begin{center}
$m_x(xy)=(m-1)+(m-1)+(m-1)^2=m^2-1$.
\end{center}
Now, the above cases are verified for $\acute{f}\in M_y(xy)$. We consider that $d(\acute{f},y)=0$ for $\acute{f}, f \in E(S(A))$. We have $m_y(xy)=1$ by $deg(y)=2$. Let $d(\acute{f},y)=1$. It yields $m_y(xy)=m-1$ and also for $d(\acute{f},y)=2$ it is obtained that $m_y(xy)=m-1$. It is shown that $m_y(xy)=2m-1$. We conclude that 
\begin{center}
$\frac{m_x(xy)}{m_y(xy)}$ = $\frac{m^2-1}{2m-1}\in \{\lambda, \frac{1}{\lambda}\}$.
\end{center}
Similarly, for $f=yz\in E(S(A))$, $m_z(yz)=m^2-1$ and $m_y(yz)=2m-1$.\\
Therefore, this results the graph $S(A)$ is $\lambda-EQDB$, so it is completed.
\end{proof}

In the following figure, we illustrate an example of the Proposition 4.2 for $A=K(3,3)$.
\begin{figure}
\begin{center}
\includegraphics[scale=0.5]{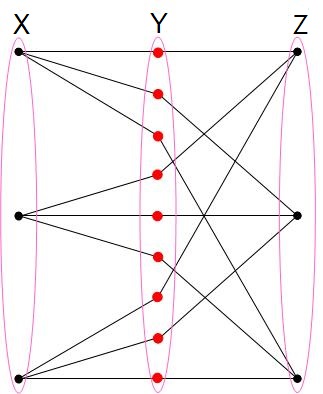}
\caption{$Graph$ $S(A)$}
\label{fig:3}
\end{center}      
\end{figure}

\begin{corollary}
Let $A$ be a regular complete bipartite graph. For any edge $f=gh\in S(A)$, it holds that
\begin{center}
$m_0(gh)=(m-1)^2$.
\end{center}
\end{corollary}

\begin{proof}
Simply it is seen $|E(S(A))|-(m_g(gh)+m_h(gh)-1)=m_0(gh)$.
\end{proof}

In this part, another exapmple of $\lambda-EQDB$ graphs is given. At the begining, the below definition is presented.
\begin{definition}
Suppose that $A$ is a connected graph. $O(A)$ is constructed from $A$ by adding two new vertices corresponding to any edge of $A$, then connecting every two new vertices to the end vertices of the corresponding edge and deleting it. Another procedure to describe $O(A)$ is to exchange each edge of $A$ with a square. (for an example see Fig 1, 4).
\end{definition}

\begin{figure}
\begin{center}
\includegraphics[scale=0.5]{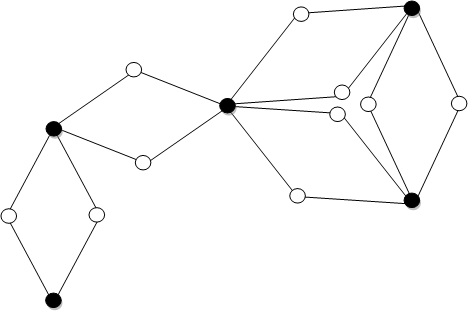}
\caption{$Graph$ $O(A)$}
\label{fig:4}   
\end{center}   
\end{figure}

In a completely similar way to the Proposition 4.2, we can conclude the following proposition and refuse to prove it. 
\begin{proposition}
Assume that $A$ is a regular complete bipartite graph of order $m$, where $m\geqslant3$. Then the graph $O(A)$ is $\lambda-EQDB$ if and only if $\lambda=\frac{2m^2-1}{4m-3}$ and also for every edge $f=gh$ we have $m_0(gh)=2m^2-4m+3$.
\end{proposition}

\section{$\lambda-NQEDB$ and $\lambda-SQEDB$ graphs}

If, for every arbitrary edge $f\in E(A)$, the constant of $m_h^A(f)$ is $\acute{\gamma}_A,$ then $A$ is introduced as \textit {Nicely edge quasi-$\lambda$-distance-balanced} ($\lambda-NQEDB$). Similarly, if for every adjacent pair of vertices $f,g$ the constant of $n_h^A(f)$ is $\gamma_A$, then $A$ is named \textit{Nicely quasi $\lambda$-distance-balanced} ($\lambda-NQDB$).\\
Each $\lambda-NQEDB$ graph is $\lambda-QEDB$, and every bipartite $\lambda-QEDB$ graph is $\lambda-NQEDB$. Let's start with the following lemma.
\begin{lemma}
If $A$ is a connected $\lambda-NQEDB$ graph with a diameter of $d$, then for every arbitrary edge $f=gh\in E(A)$, there will precisely be $|E(A)|-(\lambda+1)\acute{\gamma}_{A}$ edges of $A$ that will be at the equal distance from $g$ and $h$. On the other hand, $m_0^A(f)=|E(A)|-(\lambda+1)\acute{\gamma}_A-1$.
\end{lemma}
\begin{proof}
It is directly deduced from $m_g^A(f)+m_h^A(f)+m_0^A(f)=|E(A)|-1$.
\end{proof}

\begin{lemma}
Let $A$ be a connected $\lambda-NQEDB$ graph with a diameter $d$. Then $d-1\leq \lambda\acute{\gamma}_A$.
\end{lemma}

\begin{proof}

Consider a path $a_0,a_1,a_2,...,a_{d-1},a_d$. As well as, we have $h_1,h_2,h_3,...,h_d$ as a series of edges, where $d_A(h_1,h_d)=d$. Now, consider $h_1=a_0a_1$ is an arbitrary edge of $A$. without loss of generality, we presume that $m_{a_1}^A(h_1)=\lambda m_{a_0}^A(h_1)$. Therefore, $|\{h_2,...,h_d\}|=d-1\leq m_{a_1}^A(h_1)=\lambda m_{a_0}^A(h_1)$. It is seen that  $d-1\leq \lambda\acute{\gamma}_A$.
\end{proof}

\begin{lemma}
Assume that $A$ and $B$ are graphs. Then $A\square B$ is $\lambda-NQEDB$ if and only if both $A$ and $B$ are both $\lambda-NQEDB$ and $\lambda-NQDB$ and
\begin{center}
$|V(B)|.\acute{\gamma}_A+|E(B)|.\gamma_A=|V(A)|.\acute{\gamma}_B+|E(A)|.\gamma_B$.
\end{center}
\end{lemma}

\begin{proof}
Consider that vertices $(a_1,b_1)$ and $(a_2,b_2)$ are in $V(A\square B)$. According to the definition, either $b_1=b_2$ and $a_1$ and $a_2$ are adjacent in $A$, or $a_1=a_2$ and $b_1$ and $b_2$ are adjacent in $B$. Let now $b_1=b_2$ and $a_1$ and $a_2$ be adjacent in $A$. Assume now that $A\square B$ is $\lambda-NQEDB$. By equations (3.2),(3.3),(3.4) and (3.5), we obtain
\begin{center}
$m_{a_1}^A(a_1a_2).|V(B)|+n_{a_1}^A(a_1a_2).|E(B)|=\lambda\{m_{a_2}^A(a_1a_2).|V(B)|+n_{a_2}^A(a_1a_2).|E(B)|\}$=
$m_{b_1}^B(b_1b_2).|V(A)|+n_{b_1}^B(b_1b_2).|E(A)|=\lambda\{m_{b_2}^B(b_1b_2).|V(A)|+n_{b_2}^B(b_1b_2).|E(A)|\}$=
\end{center}

$\hspace*{4cm}$ $t(\acute{\gamma}_{A\square B}+\gamma_{A\square B})$,$\hspace*{4.30cm}(5.6)$\\
where
\begin{center}
$\gamma_{A\square B}=n_{a_1}^A(a_1a_2).|E(B)|=n_{a_2}^A(a_1a_2).|E(B)|$,
\end{center}
and
\begin{center}
$\acute{\gamma}_{A\square B}=m_{a_1}^A(a_1a_2).|V(B)|=m_{a_2}^A(a_1a_2).|V(B)|$.
\end{center}
We conclude that from the above equation that $A$ and $B$ are both $\lambda-NQEDB$ and $\lambda-NQDB$ and
$|V(B)|.\acute{\gamma}_A+|E(B)|.\gamma_A=|V(A)|.\acute{\gamma}_B+|E(A)|.\gamma_B$ satisfies.
Conversely, if $A$ and $B$ are both $\lambda-NQEDB$ and $\lambda-NQDB$ with $|V(B)|.\acute{\gamma}_A+|E(B)|.\gamma_A=|V(A)|.\acute{\gamma}_B+|E(A)|.\gamma_B$. By equations (3.2),(3.3),(3.4) and (3.5), $A\square B$ is a $\lambda-NQEDB$ graph.
\end{proof}

While for every arbitrary edge $f=gh \in E(A)$ and every $i\in[1,d-1]$ in a graph $A$ with a diameter $d$, we have $|\acute{D}_{i-1}^i(f)|=\lambda|\acute{D}_i^{i-1}(f)|+(\lambda-1)$, then such graphs are called \textit{Srtongly edge quasi $\lambda$-distance-balanced} ($\lambda-SQEDB$). Let $\lambda=1$. Then graph $A$ is a $SEDB$ graph.
\begin{example}
A class of $\lambda-SQEDB$ graphs are complete bipartite graphs $K_{n,\lambda n}$.
\end{example}

\begin{proof}
Consider $f=gh\in E(K_{n,\lambda n})$. Since $K_{n,\lambda n}$ is bipartite with a diameter 2, we obtain $\acute{D}_1^2(f)=(n-1)$ and $\acute{D}_2^1(f)=\lambda(n-1)$. Therefore, $\acute{D}_2^1(f)=\lambda\acute{D}_1^2(f)+(\lambda-1)$ and it completes the result.
\end{proof}

\begin{lemma}
Consider that $A$ is a $\lambda-SQEDB$ graph with a diameter 2. Then $A$ is a $\lambda-QEDB$ graph for $\lambda\in\mathbb{R}$ and $\lambda>1$.
\end{lemma}

\begin{proof} 
For every graph $\lambda-EQDB$, $m_\alpha^A(f)=\lambda m_\beta^A(f)$ holds and we have for every $i\in[1,d-1]$
 \begin{center}
$|\{f\}\bigcup^{d-1}_{i=1}\acute{D}_{i-1}^i(f)|=\lambda|\{f\}\bigcup^{d-1}_{i=1}\acute{D}_i^{i-1}(f)|$.
\end{center} 
 Therefore,
  \begin{center}
 $\sum^{d-1}_{i=1}|\acute{D}_{i-1}^i(f)|=\lambda\sum^{d-1}_{i=1}|\acute{D}_i^{i-1}(f)|+(\lambda-1)$.
 \end{center}
Consider now $A$ is a $\lambda-SQEDB$ graph. Then for every $f=gh\in E(A)$ and every $i\in[1,d-1]$, it is obtained
\begin{center}
$|\acute{D}_i^{i-1}(f)|=\lambda|\acute{D}_{i-1}^i(f)|+(\lambda-1)$.
\end{center}
Thus, 
\begin{center}
$\sum_{i=1}^{d-1}|\acute{D}_i^{i-1}(f)|=\lambda\sum_{i=1}^{d-1}|\acute{D}_{i-1}^i(f)|+(\lambda-1)(d-1)$.
\end{center}
If $d=2$, then we have
\begin{center}
$|\acute{D}_2^1(f)|=\lambda|\acute{D}_1^2(f)|+(\lambda-1)$.
\end{center}
The proof is completed.
\end{proof}

Notice that If $\lambda=1$, then every $SEDB$ graph is $EDB$.

 \noindent Zohreh Aliannejadi \\ Department of Mathematics, Islamic Azad University, South Tehran Branch, \\ Tehran, Iran. \\ e-mail: z\_alian@azad.ac.ir 
 \bigskip \noindent 
 \\Somayeh Shafiee Alamoti \\ Department of Mathematics, Islamic Azad University, South Tehran Branch, \\ Tehran, Iran. \\ e-mail: shafiee.s88@gmail.com 
 \bigskip \noindent 

\begin{thebibliography}{00} 

  \bibitem{1}
A. Abedi, M. Alaeiyan, A. Hujdurović, K. Kutnar, Quasi-$\lambda$-distance-balanced graphs, Discrete Appl. Math. 227 (2017) 21-28.
 \bibitem{2}
 Z. Aliannejadi, A. Gilani, M. Alaeiyan, J. Asadpour, On some properties of edge quasi-distance-balanced graphs, Journal of Mathematical Extension. 16 (2022) 1-13.
  \bibitem{3}
 Z. Aliannejadi, M. Alaeiyan, A. Gilani, Strongly Edge Distance-Balanced Graph Products, 7th International Conference on Combinatorics, Cryptography, Computer Sciences and Computing, Iran University of Science and Technology. November 16-17 (2022). 
\bibitem{4}
K. Balakrishnan, M. Changat, I. Peterin, S. Špacapan, P. Šparl, A. R. Subhamathi, Strongly distance-balanced graphs and graph products, European J. Combin. 30 (2009) 1048-1053.
\bibitem{5}
D.M. Cvetkoci\'{c}, M. Doob, H. Sachs, Spectra of Graphs-Theory and Application, Academic Press, New York, (1980).
\bibitem{6}
M. Faghani, A.R. Ashrafi, Revised and edge revised Szeged indices of graphs, Ars Math. Contemp. 7 (2014) 153-160. 
\bibitem{7}
A. Graovac, M. Juvan, M. Petkovsek, A. Vesel, J. Zerovnik, The Szeged index of fasciagraphs, Match Common. Chem. 49 (2003) 47-66. 
\bibitem{8}
I. Gutman, L. Popovic, P.V. Khadikar, S. Karmarkar, S. Joshi, M. Mandloi, Relations between Wiener and Szeged indices of monocyclic molecules, Match Common. Math. Comput. Chem. 35 (1997) 91-103. 
\bibitem{9}
I. Gutman,A.R. Ashrafi, The edge version of the Szeged index, Croat. Chem. Acta. 81 (2008) 263-266. 
\bibitem{10}
A. Hujdurovi\'{c}, On some properties of quasi-distance-balanced graphs, Bull. Aust. Math. Soc. 97 (2018) 177-184. 
\bibitem{11}
A. Ili\v{c}, S. Klav\v{z}ar, M. Milanovi\'{c}, On distance-balanced graphs, European J. Combin. 31 (2010) 733-737. 
\bibitem{12}
W. Imrich, S. Klav\v{z}ar, Product Graphs: Structure and Recognition, Wiley, New York, USA. (2000).
\bibitem{13}
J. Jerebic, S. Klav\v{z}ar, D.F. Rall, Distance-balanced graphs, Ann. Comb. 12 (2008) 71-79. 
\bibitem{14}
M.H. Khalifeh, H. Yousefi-Azari, A.R. Ashrafi, S.G. Wagner, Some new results on distance-based graph invariants, European J. Combin. 30 (2009) 1149-1163. 
\bibitem{15}
K. Kutnar, A. Malni\v{c}, D. Maru\v{s}i\v{c}, \v{S}. Miklavi\v{c}, Distance-balanced graphs: Symmetry conditions, Discrete Math. 306 (2006) 1881-1894. 
\bibitem{16}
K. Kutnar, A. Malni\v{c}, D. Maru\v{s}i\v{c}, \v{S}. Miklavi\v{c}, The strongly distance-balanced property of the generalized Petersen graphs, Ars Math. Contemp. 2 (2009) 41-47. 
\bibitem{17}
K. Kutnar, \v{S}. Miklavi\v{c}, Nicely distance-balanced graphs, European J. Combin. 39 (2014) 57-67. 
\bibitem{18}
\v{S}. Miklavi\v{c}, P. \v{S}parl, On the connectivity of bipartite distance-balanced graphs, European J. Combin. 33 (2012) 237-247. 
\bibitem{19}
M. Tavakoli, H. Yousefi-azari, A.R. Ashrafi, Note on Edge Distance-Balanced Graphs, Transaction on Combinatorics, University of Isfehan. 1 no. 1 (2012) 1-6. 
\bibitem{20}
W. Yan, B.Y. Yang, Y.N. Yeh, The behavior of Wiener indices and polynomials of graphs under five graph decoration, Appl. Math. Lett. 20 (2007) 290-295. 

 \end{thebibliography}
 \end{document}